\providecommand{\href}[2]{#2}

\documentclass[12pt]{amsart}

\usepackage{graphicx}
\usepackage{amsthm}
\usepackage{amsmath, amssymb}
\usepackage{thmtools}
\usepackage{thm-restate}

\usepackage{hyperref}

\usepackage{cleveref}

\swapnumbers

\theoremstyle{plain}
\newtheorem{Thm}{Theorem}

\newtheorem{Coro}[Thm]{Corollary}

\newtheorem{Claim}[Thm]{Claim}

\newtheorem{Rem}[Thm]{Remark}

\theoremstyle{definition}
\newtheorem{Def}[Thm]{Definition}

\begin{document}
\begin{abstract}  Let $K$ be a fibered knot in $S^3$. We show that if the monodromy of $K$ is sufficiently complicated, then Dehn surgery on $K$ cannot yield a lens space.  Work of Yi Ni shows that if $K$ has a lens space surgery then it is fibered.  Combining this with our result we see that if $K$ has a lens space surgery then it is fibered and the monodromy is relatively simple.     

\end{abstract}
\title{Dehn surgery on complicated fibered knots in the 3-sphere}

\maketitle

 \author Abigail Thompson \footnote{Supported by the National Science Foundation and Friends of the Institute for Advance Study.}

\section{Introduction}

Let $K$ be a knot in $S^3$.   One can obtain a new manifold $M$ by removing an open  neighborhood of $K$ and attaching a solid torus $T$ to the resulting knot complement via some homeomorphism $h$ from $\partial{T}$ to $\partial{N(K)}$.       The homeomorphism $h$ is completely determined by a pair of relatively prime integers $(p,q)$,  where $h$ maps the boundary of a meridian disk of $T$ to a curve $\alpha$ that wraps $p$ times around $K$ meridionally and $q$ times longitudinally.  This is called {\it{surgery}} on $K$.    If $q=1$ the surgery is $integral$.   The definition extends in an obvious way to a link with more than one component.

It is a useful and classical result that every 3-manifold can be obtained by  surgery on a link in $S^3$ \cite{L}, which in general has many components.  A natural question thus arises:  which manifolds can be obtained by surgery on a knot in $S^3$, and more specifically, on which knots?

If one starts with the trivial knot $K$, it is clear that any surgery yields $S^3$, $S^1\times{S^2}$ or a lens space $L$.     A restricted form of the question asks, which knots in the 3-sphere can yield one of these simple manifolds after non-trivial surgery?   

Here is a brief synopsis of some of the progress so far on this question:

In \cite{G}, Gabai showed that if surgery on $K$ yields $S^1\times{S^2}$ then the knot is the unknot (and the surgery is $(0,1)$).   In \cite{G-L}, Gordon and Luecke showed that if (non-trivial) surgery on a knot yields $S^3$, then the knot is trivial.     So the question has a complete answer for $S^3$ and $S^1\times{S^2}$.    

Berge \cite{B} produced examples of non-trivial knots $K$ such that surgery on these knots yields a lens space, and Gordon \cite{K} conjectured that his list contains all non-trivial knots with a lens space surgery.     These are now called {\it Berge knots}.  While there has been substantial progress on Berge's conjecture, it remains open.  

A major step towards a resolution of the Berge conjecture was provided by the cyclic surgery theorem \cite{C-G-L-S}, which shows that a lens space surgery must be integral.    Another major development was provided by Yi Ni \cite{YN}, who showed that a knot with a lens space surgery must be fibered.    
Obtaining $L(2,1)$, and so $(2,1)$ surgery,  was ruled out for non-trivial knots by Kronheimer,  Mrowka,  Ozsvath and Szabo \cite{K-M-O-S}.    Recently Greene \cite{Gr} has shown that if a knot has $K$ has a lens space surgery, then surgery on some Berge knot yields the same lens space.   

Compiling a few of these results, we see that if $K$ is knot in the 3-sphere such that $(p,q)$ surgery on $K$ yields a lens space, then $q=1$,  $p\geq3$ and $K$ is fibered.   We show that if the fibration of $K$ is sufficiently complicated, then $K$ does not have a lens space surgery.       

\begin{figure}[h]
    \centering
    \includegraphics[width=0.6\textwidth]{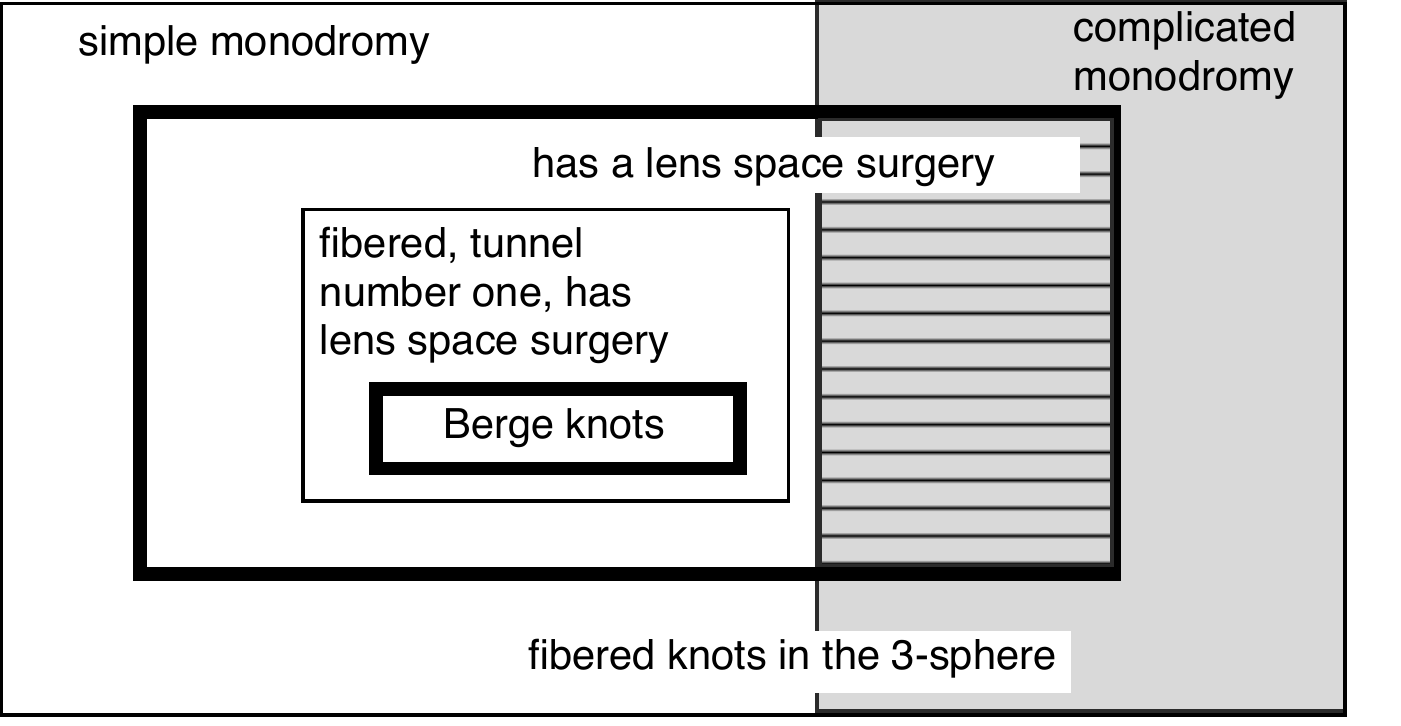}
    \caption{Knots and the Berge conjecture}
    \label{fibered3}
\end{figure}

 It follows from results of Banks and Rathbun \cite{B-R}  that the knots under consideration, namely  fibered knots $K$ in the 3-sphere with a complicated monodromy, must have tunnel number greater than one, and hence are not  Berge knots.      Thus our result is additional evidence in favor of a positive answer to the Berge conjecture.  Figure \ref{fibered3} illustrates the relations among some of these results.   Berge's conjecture amounts to showing that the two bold boxes are the same.     Our contribution is to show that the doubly-shaded region is empty.  We note that S. Schleimer conjectures \cite{Sch} the the entire shaded region is empty, that is, that all fibered knots in the 3-sphere are simple.    

An outline of the paper is as follows:

In section 2, we  prove the two easier cases of the main theorem.    These warm-up cases show that $(p,1)$ ($p\geq3$) surgery on a knot $K$ in the 3-sphere cannot yield either a reducible manifold or a toroidal manifold.    The strategy is to intersect a 2-complex constructed from the boundary of $K$ and two copies of the fiber with the supposed reducing 2-sphere or essential torus, and, using the complication of the monodromy, to derive a contradiction from the graph of intersection.     

In section 3 we prove the final case, in which we show that $(p,1)$ ($p\geq3$) surgery on a knot $K$  cannot yield a lens space.   The strategy is as before, but more delicacy is required.   For the results in section 2, we use a 1-parameter sweep-out of the fiber surface relative to the reducing 2-sphere or essential torus.     For the proof of the final case, we need  two 1-parameter arguments, using both the fibration of the knot complement and the sweep-out induced by the supposed genus one Heegaard splitting of the resulting manifold to arrive at a similar contradiction.   In Section 4 we state some extensions to knots in arbitrary 3-manifolds suggested by David Futer.

The author is grateful to Marty Scharlemann for his extensive help with this paper.  David Futer pointed out that these results are easily extended to fibered knots in arbitrary 3-manifolds, and to non-integral surgeries on knots in the 3-sphere.  The restriction on the surgery curve is simply that it must intersect the boundary of the fiber in at least three points.   Finally  the author is grateful to the Institute for Advanced Study, where she was the  Neil Chriss and Natasha Herron Chriss Founders' Circle Member for 2015-2016 while this work was being completed.     

\section{Warm-up cases}

Let $K$ be a non-trivial fibered knot in $S^3$ with fiber $F$ over $S^1$.    We parameterize $S^1$ by $(\cos{2\pi{t}},\sin{2\pi{t}}) $ with $t\in[0,1]$.   Let $F_0$ and $F_t$  be two disjoint copies of the fiber, corresponding to inverse images of  $0$ and $t\neq{0}$.      $F_0$ and $F_t$ split $K^c$ into two handlebodies, $A_t$ and $B_t$, each of which has a natural $F\times{I}$  structure.      Let $h$ be the monodromy of the fibration.      

\begin{Def} A product disk for $A_t$ ($B_t$) is a properly imbedded disk which is the product of an essential arc $\alpha$ properly imbedded in $F_0$ ($F_t$) with $I$.    A product annulus for $A_t$ ($B_t$) is an properly imbedded annulus in $A_t$ ($B_t$) which is the product of an essential simple closed curve $\gamma$ imbedded in $F_0$ ($F_t$) with $I$.  We say that the monodromy $h$ is {\it complicated} if any pair $(Q_A,Q_B)$ of product disks/annuli in $(A_t,B_t)$ intersect.  Otherwise we say the monodromy is {\it simple}. 
\end{Def}

Being simple is related to the monodromy having low translation distance in the arc complex of the fiber (see \cite{S-Y} \cite{B-R} \cite{F-S} for definitions), but the definitions of translation distance in this context require various degrees of  control on the endpoints of arcs.

{\it Simple} most closely corresponds to {\it translation distance at most two in the arc and curve complex} (see \cite{F-S}, \cite{H-P-W}).   To emphasize the geometric aspect of splitting the knot complement into two handlebodies, we will stick with the simple/complicated terminology.   

\begin{Thm}\label{Main}   Let $K$ be a fibered knot in $S^3$ with monodromy $h$.   If $h$ is complicated, then p-Dehn surgery on $K$, $p\geq3$,  does not yield (1) a reducible manifold, (2) a toroidal manifold, or (3) a lens space.

\end{Thm}

Combining this with previous results (\cite{C-G-L-S}, \cite{YN},\cite{K-M-O-S}) gives the immediate corollary:

\begin{Coro} Let $K$ be a knot in $S^3$.   If surgery on $K$ yields a lens space, then $K$ is fibered and the monodromy of $h$ is simple.
\end{Coro}

Proof of Theorem:

Case 1:
Assume $p\geq3$  surgery on $K$ yields a reducible manifold $M$.   Let $K'$ be the dual of $K$ in $M$, i.e., $K'$ is the core of the surgered solid torus.   Let $S$ be a reducing sphere for $M$. In both case 1 and 2, we use the fact that $F=F_0$, $S$ and $T$ are essential in $K^c$ to  remove non-essential simple closed curves of intersection between $F$ and $S$ (Case 1) or $T$ (Case 2), while making $K'$ intersect $S$ ($T$) transversely.   $K'\cap{S}\neq\emptyset$ since $S^3$ is not reducible; $K'\cap{T}\neq\emptyset$ since $h$ is not simple (if $T$ lies in the complement of $K'$ and hence of $K$, then it can be isotoped to intersect $F$ in essential simple closed curves, which define annuli contradicting that $h$ is simple) .    

Imbedded in $M$ is a family of two complexes $W_t$, where $W_t=(\partial{N(K)\cup{F_0}\cup{F_t}})$.    We consider the trivalent graph  $\Gamma_t$ in $S$ obtained by intersecting $S$ with $W_t$.     $W_t$ splits $M$ into three handlebodies,   $A_t$, which we will color black, $B_t$, which we color white, and the solid torus bounding the surgery curve, which we color red.     Thus the regions of  $\Gamma_t$ in $S$ are similarly colored black, white and red.     

\begin{figure}[h]
    \centering
    \includegraphics[width=0.4\textwidth]{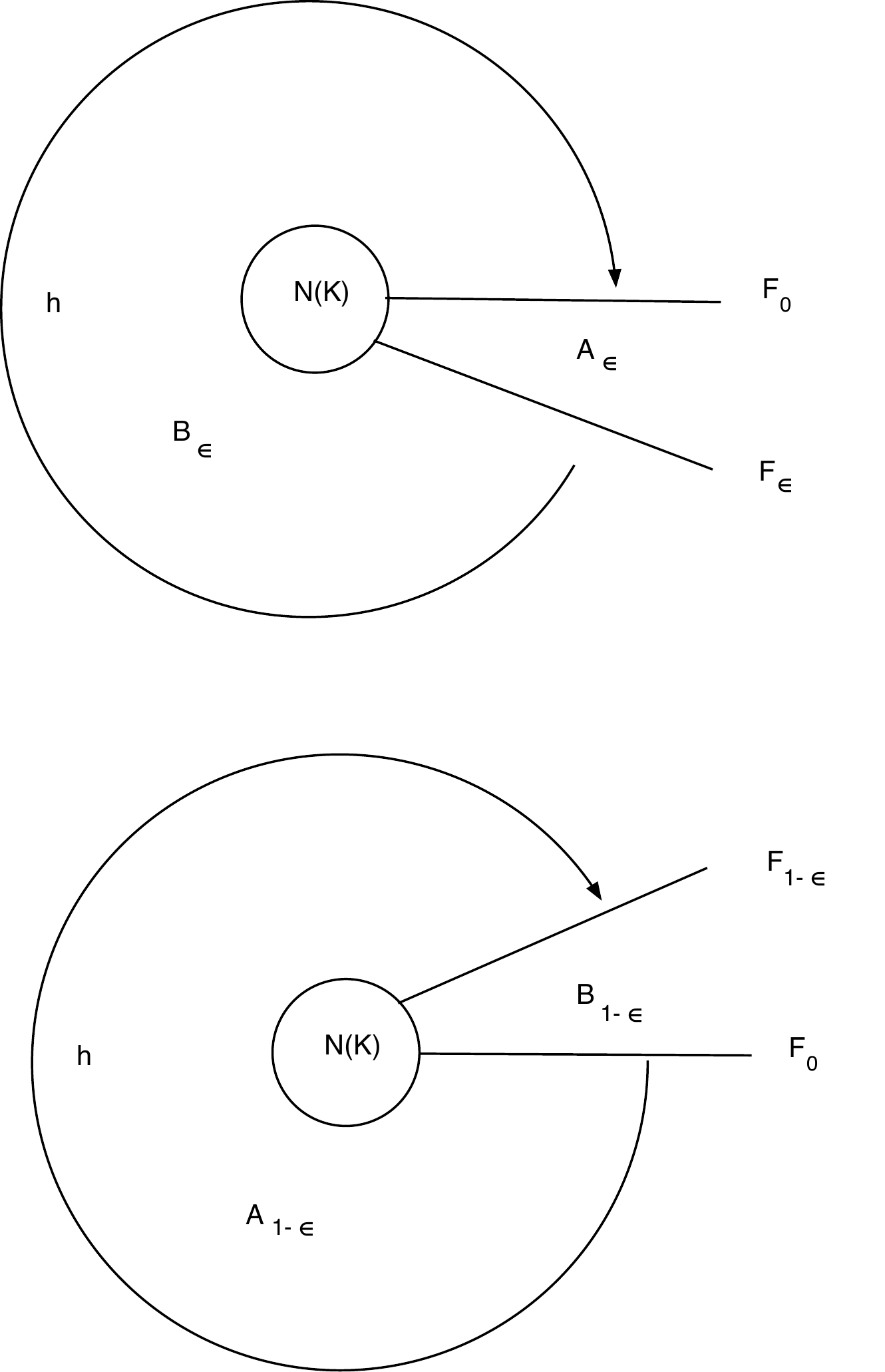}
    \caption{$W_{\epsilon}$ and $W_{1-\epsilon}$}
    \label{fibered2}
\end{figure}

\begin{figure}[h]
    \centering
    \includegraphics[width=0.8\textwidth]{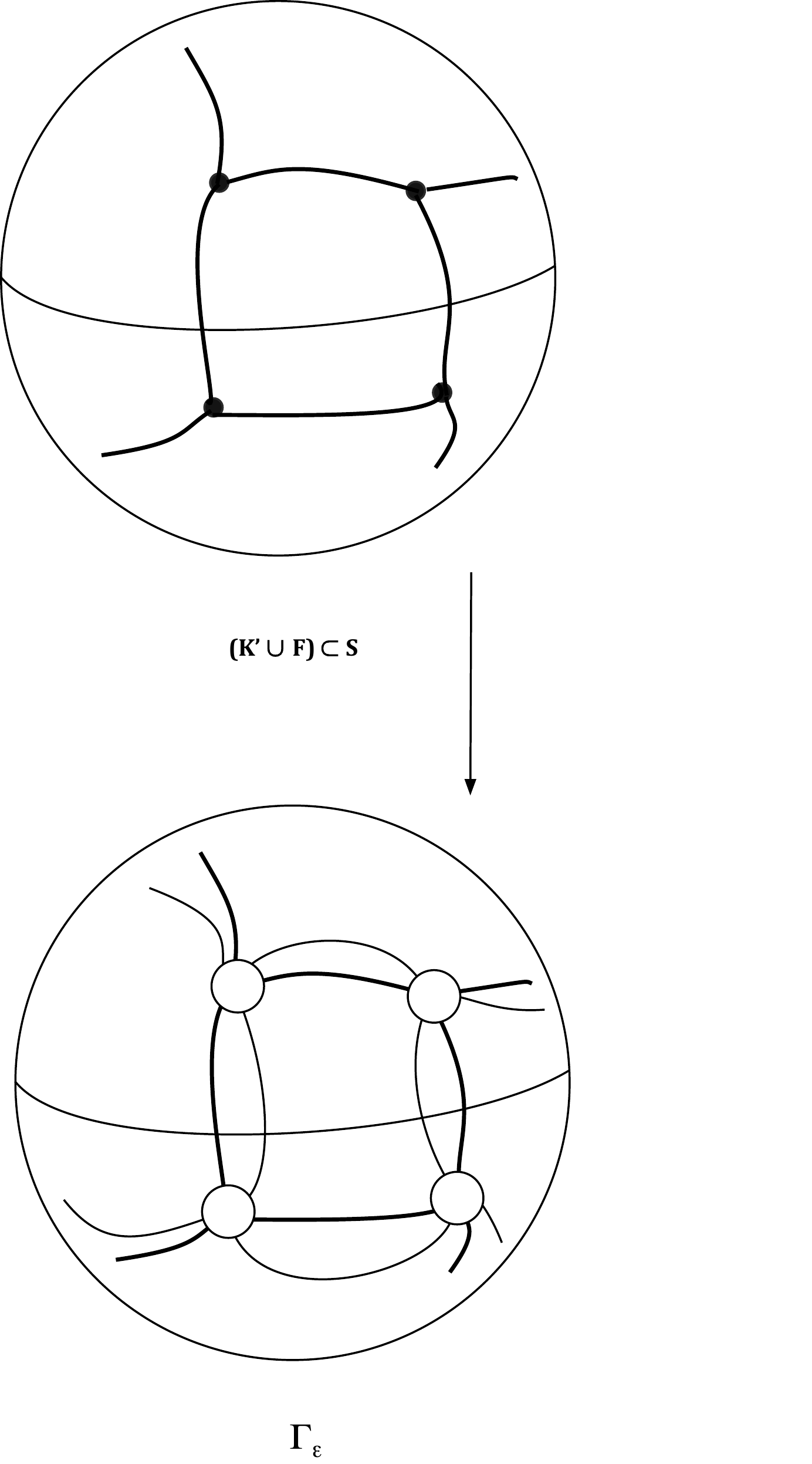}
    \caption{constructing $\Gamma_{\epsilon}$}
    \label{fibered4}
\end{figure}

Notice that for $\epsilon$ very small, the black regions of  $\Gamma_{\epsilon}$ in $S$ are all product disks or annuli in $A_\epsilon$.     Similarly, the white regions of $\Gamma_{1-\epsilon}$ in $S$ are all product disks or annuli in $B_ {1-\epsilon}$ (see Figures \ref{fibered2} and \ref{fibered4}).   

Since $h$ is complicated, we never see white and black product regions either appearing simultaneously in $\Gamma_t$ or ``switching off'' between $\Gamma_t$ and $\Gamma_{t+\epsilon}$ (see \cite{G}, Lemma 4.4, and Figure \ref{fibered11}).    Thus for some value $t_1$ of $t$, no region of  $\Gamma_{t_1}$ is a product disk or annulus.     

If $\Gamma_{t_1}$ contains a simple closed curve or is disconnected, we pass to an innermost component; for simplicity we will continue to refer to this as  $\Gamma_{t_1}$.

Let $V$, $E$ and $F$ be the number of vertices, edges and faces  in $\Gamma_{t_1}$.   Since $\Gamma_{t_1}$ is trivalent, $E={\frac{3}2}V$.    Since there are no white or black product regions, each white/black region has at least $6$ edges.    Since the surgery on $K$ is at least $3$, each red region has at least $6$ edges, so $E\geq{3F}$.    An Euler characteristic calculation on the disk shows this is impossible, since $V-E+F\leq{\frac{2}3}E-E+{\frac{1}3}E=0$.   Hence there is no such reducing sphere $S$. 

Case 2:
Assume $p\geq3$  surgery on $K$ yields a toroidal manifold $M$.   Let $T$ be an essential torus in $M$.    We define $\Gamma_{t_1}$ as before, by intersecting the torus $T$ with the  2-complex $W_{t_1}$, selecting ${t_1}$ so that no region of $\Gamma_{t_1}$ in $T$ is a product disk or annulus.    As before, if there are simple closed curves of intersection we pass to a region (now possibly an annulus) which contains no such curves.     We use the same Euler characteristic argument as in case 1, but note that we have a little more information we have not yet used.    Namely, we consider a white (black) region $X$ in $\Gamma_{t_1}$.   By the choice of ${t_1}$, we know this is not a product region, so it has more than four sides.   The singular set of $W_t$ consists of two longitudes on $\partial{N(K)}$.    We can arrange  that the curves of intersection between $T$ and $\partial{N(K)}$ are a collection of simple closed curves intersecting these longitudes both algebraically and geometrically $p$ times.    Then some sides of the disk $X$ lie on $\partial{N(K)}$, some lie on $F_0$ and some lie on $F_{t_1}$.    If we read around the boundary of $X$, we will see a side on $F_0$, followed by a side on $\partial{N(K)}$, followed by a side on $F_{t_1}$, and so on.   Thus if the boundary of $X$ has more than four sides, it must have at least eight.        Then every red region has exactly $2p$ sides, and every white or black region has at least eight sides.     Let $F_r$ be the number of red regions, $F_b$ the number of black regions and $F_w$ the number of white regions in $\Gamma_{t_1}$.    Let $V$, $E$ and $F$ be the number of vertices, edges and faces  in $\Gamma_{t_1}$.   As before,  $E={\frac{3}2}V$.     However with this more careful count on the faces, we see that $E\geq3F_r+4F_w+4F_b>3F$.   So 
$$
0=V-E+F<{\frac{2}3}E-E+{\frac{1}3}E=0,
$$

a contradiction.   Hence there is no such essential torus $T$.

\section{Main theorem, Case 3}
Assume $p\geq3$  surgery on $K$ yields a lens space $M$.   Let $T$ be a Heegaard torus in $M$.    We wish to define $\Gamma_{t_1}$ as before, by intersecting the torus $T$ with the  2-complex $W_{t_1}$, and selecting ${t_1}$ so that no region of $\Gamma_{t_1}$ in $T$ is a product disk or annulus, but in this case we have to choose $T$ with more care.     

Let $C_0$ and $C_1$ be cores of the (two) solid tori defined by the Heegaard torus $T$.   The complement of $C_0\cup{C_1}$ can be foliated by copies of $T$,   parameterized by $0<{r}<1$.      We define a function $j:M\rightarrow{[0,1]}$  induced by this foliation as follows: $j^{-1}(0)=C_0$ and $j^{-1}(1)=C_1$ , and  $j^{-1}(r)=T_r$ is a leaf of the foliation, isotopic to $T=T_{\frac{1}{2}}$, for  $0<r<1$.

As before, let $K'$ be the dual of $K$ in $M$, i.e., $K'$ is the core of the surgered solid torus.    We will find our desired $T$ in two 1-parameter steps.

Step 1:  First place $K'$ in thin position with respect to the function $j$.   For details on this procedure, see for example   \cite{H-R-S}.    Let $\Delta_r$ be the graph of intersection between $T_r$ and $F$($=F_0)$ contained in $T_r$.  

\noindent Using standard thin position arguments, we can find a level $r$ such that either  
 
\begin{itemize}
\item $K'$ lies in $T_r$ or \\
\item $K'$ intersects $T_r$ transversely, every arc in $\Delta_r$ is essential in $F$ and every simple closed curve in $\Delta_r$ is either essential both in $F$ and in $T_r$ or is inessential in both.    \\  

\end{itemize}

If $K'$ lies in $T_r$ , then $T_r\cap{K^c}$ is an essential annulus $A$ properly imbedded in the complement of $K$.   By standard innermost disk arguments, $A$ can be isotoped to intersect $F$ in essential arcs, thus cutting $A$ into disks which can be used to contradict the simplicity of $h$.    So assume we can find a level $T_r$ so that  $K'$ intersects $T_r$ transversely,  every arc in $\Delta_r$ is essential in $F$ and every simple closed curve in $\Delta_r$ is either essential both in $F$ and in $T_r$ or is inessential in both.    Ignore the latter set of simple closed curves.    

Step 2:   Using this $T_r$ (which we will now call $T$), we apply the reasoning of Case 2.    First reconstitute the 2-complex $W_t$ by taking the boundary of a (small) neighborhood of $K'$ and attaching $F$ and $F_{\epsilon}$.    This exactly transforms the graph  $\Delta_r$ into a graph $\Gamma_{\epsilon}$ contained in $T$.   As before, we examine $\Gamma_t$ as $t$ increases from $0$ to $1$.      Since the arcs and circles of $F=F_0$ are all essential ( because of Step 1), and they are unchanged as $t$ increases, any rectangles or annuli that appear in $\Gamma_t$ are necessarily product disks or annuli.      As before,  there must exist a value $t_1$ so that $\Gamma_{t_1}$ contains no product disks or annuli, and we obtain the same contradiction, using Euler characteristic, as in Case 2.     Hence there is no Heegaard torus $T$.

\section{An alternate argument, and fibered knots in arbitrary 3-manifolds}

In this section we present an alternate argument, which yields additional, somewhat more technical,  results.

\begin{Rem}\label{remark}

We can observe from the previous cases that one can obtain the required pair $(D_1,D_2)$ of disjoint product disks  from the graph $\Gamma_0$ provided just two conditions are met:\\

\noindent 1. the arcs of $\Gamma_0$ are essential in $F$.\\
\noindent 2. $\Gamma_0$ contains a disk component $D$.    
\end{Rem}

\smallskip

\begin{figure}[h]
    \centering
    \includegraphics[width=0.6\textwidth]{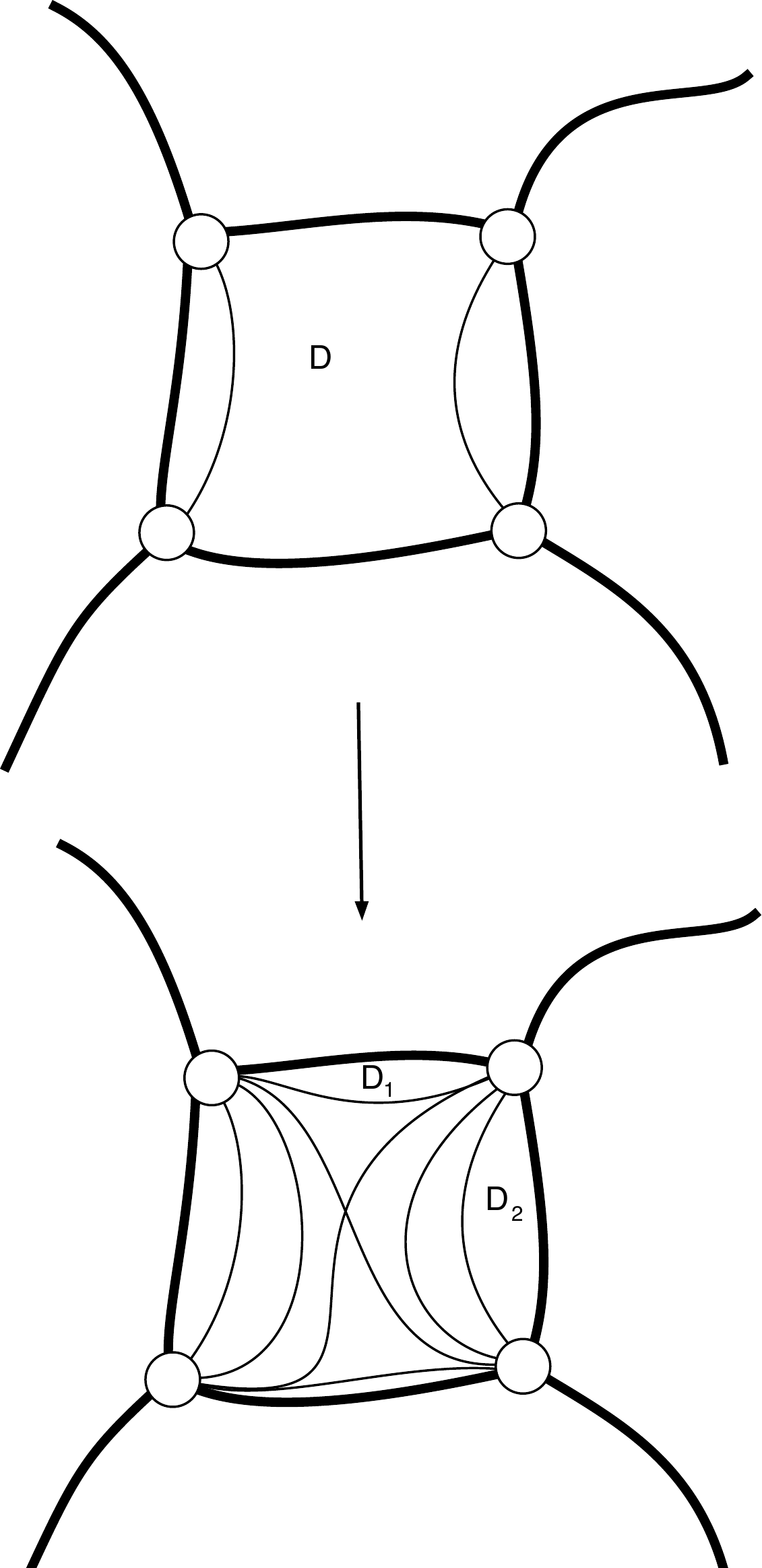}
    \caption{Singular foliation of D}
    \label{fibered6}
\end{figure}

\begin{figure}[h]
    \centering
    \includegraphics[width=0.6\textwidth]{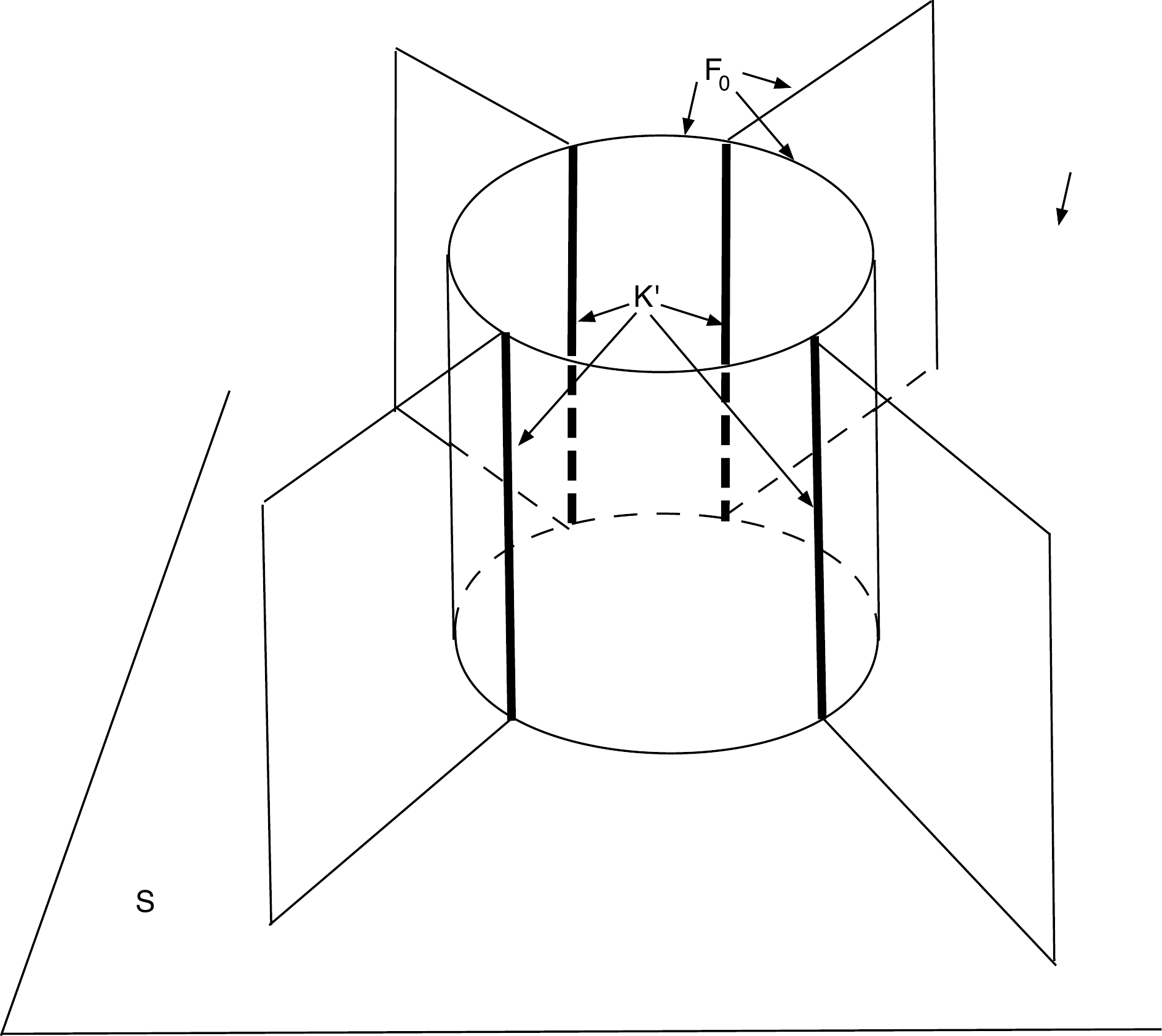}
    \caption{$K'\cup{F_0}$}
    \label{fibered10}
\end{figure}

\begin{figure}[h]
    \centering
    \includegraphics[width=0.6\textwidth]{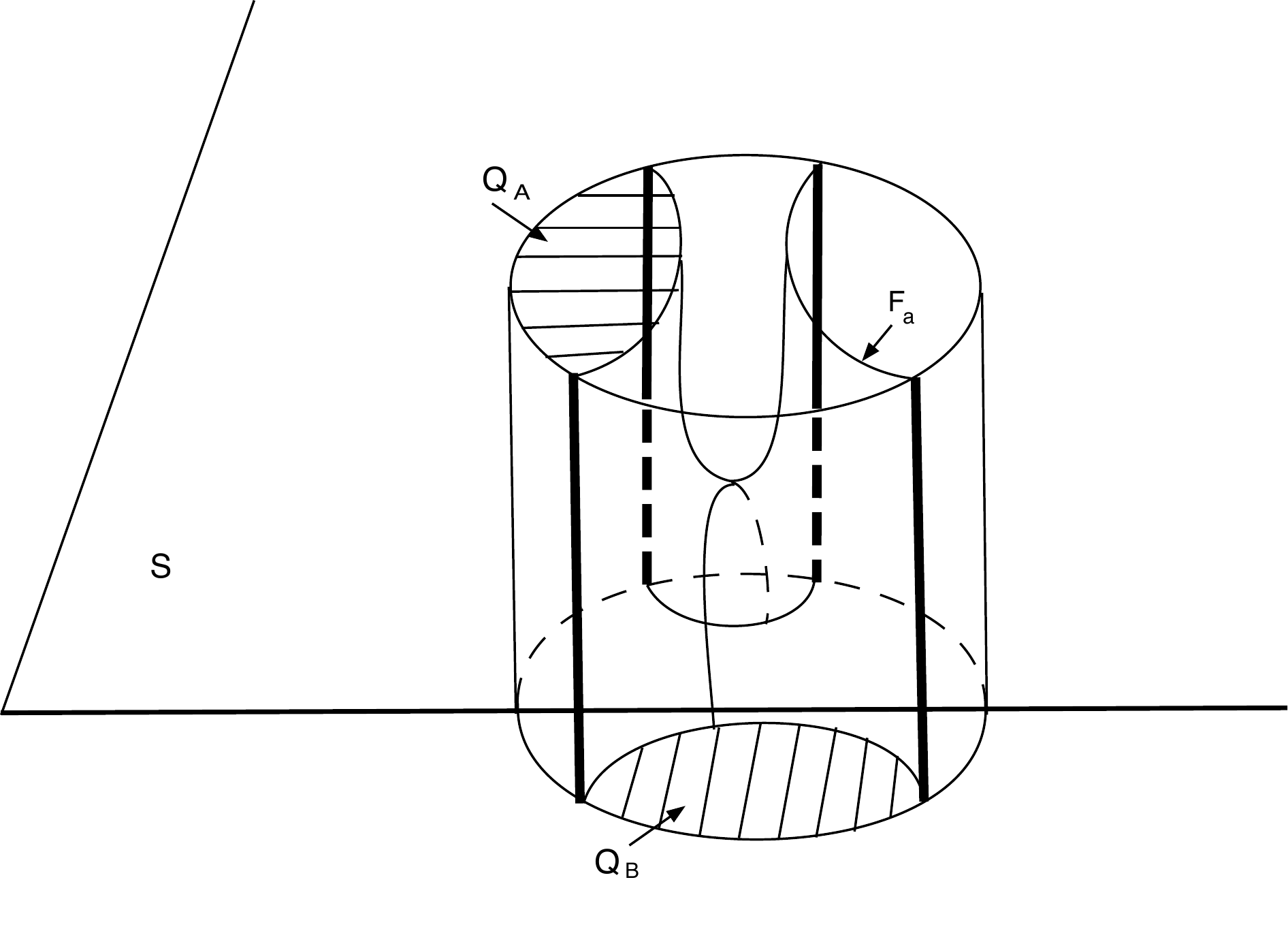}
    \caption{Disjoint product disks at a saddle exchange}
    \label{fibered11}
\end{figure}

\noindent Explanation of Remark \ref{remark}:

The fibration induces a singular foliation of $D$, by intersecting $D$ with the surfaces $F_t$.     Each non-singular level of this foliation has at least two outermost arcs, which represent product disks between $F_t$ and $F_0$ on one side of $F_t$ or the other.   When $t$ is close to zero, the visible product disks lie in $A_t$.   When $t$ is close to $1$, the visible product disks lie in $B_t$.    Hence there is either a non-singular level at which product disks on both sides appear simultaneously (an obvious pair of disjoint product disks), or there is a singular level $F_{p}$ at which a product disk on one side is exchanged through a saddle move to a product disk on the other.   This pair of product disks also constitutes a disjoint pair, as required.    See Figures \ref{fibered6},\ref{fibered10},\ref{fibered11}.

This yields the statement of two additional theorems, which, taken together, are stronger than the main theorem:

\begin{Thm}\label{Second}   Let $K$ be a fibered knot in $M^3$ with monodromy $h$.   Assume $h$ is complicated. Let $\beta$ be the slope of the boundary of the fiber on $\partial(K)$.   Let $N^3$ be obtained by $\gamma$  Dehn surgery on $K$, and let $d$ be the algebraic intersection number of $\gamma$ with $\beta$.   Let $G$ be an incompressible (closed, orientable) surface  in $N^3$.   Let $K'\subset{N^3}$ be the dual knot of $K$.   Let $g=genus(G)$.
Then at least one of the following holds:

1.  $K'$ can be isotoped to be disjoint from $G$ (notice this means $G$ is an incompressible surface in the complement of $K$ itself).

2. $d<2g+2$.

\end{Thm}

\begin{Thm}\label{Third}   Let $K$ be a fibered knot in $M^3$ with monodromy $h$.   Assume $h$ is complicated. Let $\beta$ be the slope of the boundary of the fiber on $\partial(K)$.   Let $N^3$ be obtained by $\gamma$  Dehn surgery on $K$, and let $d$ be the algebraic intersection number of $\gamma$ with $\beta$.   Let $G$ be a strongly irreducible Heegaard surface  in $N^3$.   Let $K'\subset{N^3}$ be the dual knot of $K$.   Let $g=genus(G)$.

Then at least one of the following holds:

1.  $K'$ can be isotoped to lie on $G$.

2. $d<2g+2$.

\end{Thm}

The proof proceeds as in the main theorem, using the observation above and the following claim:

\begin{Claim}   Let $\Gamma$ be a $d$-regular graph in a closed orientable surface $G$ of genus $g$.   Assume (for convenience) that $\Gamma$ has an even number of vertices;   let $2p>0$ be the number of vertices of $\Gamma$.    If $(pd)-(2p-1)\geq{2g}$, then some component of the complement of $\Gamma$ is a disk.      
\end{Claim}

Proof of Claim:  This follows by induction on the genus of $G$, using Euler characteristic.    It quantifies the notion that a graph on a surface of genus $g$ with ``enough" edges must have at least one complementary region a disk.      A more-crude estimate based on the claim then is:

\begin{Coro}\label{fourth}

With $\Gamma$ as above, if $d\geq{2g+2}$, then some component of the complement of $\Gamma$ is a disk.  

\end{Coro}

Proof of Theorems:

  $K'$ can either be isotoped disjoint from $G$ (case 1 of Theorem \ref{Second}) or to lie on $G$ (case 1 of Theorem \ref{Third}) or we can arrange $K'$ and $F$ so that arcs of the graph of intersection $\Gamma_0$ are essential in $F$.  If $d\geq{2g+2}$, Corollary \ref{fourth}  shows that this graph must contain a disk component.    Remark \ref{remark} completes the argument.

\begin{Coro}\label{Fifith}   Let $K$ be a hyperbolic fibered knot in $M^3$ with monodromy $h$. Assume $h$ is complicated.  Let $\beta$ be the slope of the boundary of the fiber on $\partial(K)$.   Let $N^3$ be obtained by $\gamma$  Dehn surgery on $K$, and let $d$ be the algebraic intersection number of $\gamma$ with $\beta$.  Let $K'$ be the dual knot of $K$ in $N^3$.

At least one of the following holds:\\
\noindent 1. $N^3$ is hyperbolic. \\
\noindent 2. $d<{6}$.\\
\noindent 3. $N^3$ is a genus 2 small Seifert fibered space and $K'$ lies on a minimal genus Heegaard surface for $N^3$.

\end{Coro}

Proof of Corollary:

If $d\geq{6}$, then $N^3$ is irreducible, atoroidal, and not a lens space, using the fact that $h$ is complicated and applying Theorems \ref{Second} and \ref{Third}.    If  $d\geq{6}$ and $N^3$ is a small Seifert fibered space, Theorem \ref{Third} implies that $K'$ can be isotoped to lie on a minimal genus Heegaard surface for $N^3$, which is possibility 3.

\begin{flushright}
Abigail Thompson\\ 
Department of Mathematics\\
University of California\\ Davis,
CA 95616\\ e-mail: thompson@math.ucdavis.edu\\
\end{flushright}


\begin{thebibliography}{HHH}


\bibitem{B-R} J.E. Banks, M. Rathbun, {\it Monodromy action on unknotting tunnels in fiber surfaces},  arXiv:1312.6902v2 [math.GT] .


\bibitem{B} J. Berge, {\it Some knots with surgeries yielding lens spaces},  Unpublished manuscript, c. 1990.

\bibitem{C-G-L-S} M. Culler, C. McA. Gordon, J. Luecke, and P. B. Shalen, {\it Dehn surgery on knots}, Bull. Amer. Math. Soc. (N.S.) Volume 13, Number 1, 43-45 (1985).

\bibitem{F-S} D. Futer, S. Schleimer, {Cusp geometry of fibered 3-manifolds}, Amer. J. Math. 136, no. 2, 309-356, (2014).

\bibitem{G} D. Gabai, {\it Foliations and the topology of 3-manifolds. III}, J. Differential Geom.
Volume 26,  479-536 (1987).





\bibitem{G-L} C. McA. Gordon  and J. Luecke, {\it Knots are determined by their complements}, J.A.M.S., vol 2.,  number 2, 371-415, (1989).

\bibitem{Gr} J. Greene, {\it The lens space realization problem}, Annals of Mathematics 177, 449-511, (2013).

\bibitem{H-P-W} S. Hensel, P. Przyticki, R. C. H. Webb, {\it 1-slim triangles and uniform hyperbolicity for arc graphs and curve graphs}, arXiv:1301.5577.  Journal of the European Mathematical Society, vol. 17, number 4, 755?762 (2015).
 
 
\bibitem{H-R-S} H. Howards, Y. Rieck, J. Schultens, 

{\it Thin position for knots and 3-manifolds:  a unified approach}, Geometry and Topology Monographs 12,  89-120 (2007). 

\bibitem{K} R. Kirby, {\it Problems in low-dimensional topology}, math.berkeley.edu/~kirby/problem 1.78, (2010).

\bibitem{K-M-O-S} P. Kronheimer, T. Mrowka, P. Ozsvath and Z. Szabo, {\it Monopoles and lens space surgeries},  Ann. of Math. 165(2):457-546, (2007).





\bibitem{L} W. B. R. Lickorish, {\it A representation of orientable combinatorial 3-manifolds}, Ann. of Math. 76 (3): 531-540,  (1962).

\bibitem{YN}  Y. Ni, {\it Knot Floer homology detects fibred knots}, Inventiones Mathematicae 170 (3): 577-608, (2007).

\bibitem{S-Y} T. Saito, R. Yamamoto, {\it Complexity of open book decompositions via the arc complex}, Journal of Knot Theory and Its Ramifications 19 (01): 55-69,(2010).

\bibitem{Sch} S. Schleimer, {\it private communication}.   


\bibitem{W} A. H. Wallace, {\it Modifications and cobounding manifolds}, Canad. J. Math. 12: 503-528, (1960).

%
 \end{thebibliography}
\end{document}